%I. Please use the following macros in your TeX file.

\input amssym
\magnification=\magstephalf
\voffset=30pt
\hoffset=55pt
\parskip=.5pt
\parindent=0pt
\hsize=5.00 true in
\vsize=7.4 true in
\font\smallit=cmti9

%\headline={\ifnum\pageno<2 \hfil \else \smalltt INTEGERS: 21 (2021), \#Axx \hfil \folio \fi} % We will fill in the x's

\footline={\hfil\rm\folio\hfil}

%II. For the title page (title, authors, and abstract) please use the following %construct.

\centerline{\bf CONGRUENCES FOR SUMS OF POWERS OF AN INTEGER}
\vskip 30pt
\centerline{{\bf Leif Jacob}%\footnote{$^1$}{any footnote here}
}
\centerline{\smallit Institute for Mathematics, Friedrich Schiller University,  
Jena, Germany}
\centerline{\tt leif.jacob@uni-jena.de} %(optional)
\vskip 10pt
\centerline{{\bf Burkhard K\"ulshammer}%\footnote{$^2$}{any footnote here}
}
\centerline{\smallit Institute for Mathematics, Friedrich Schiller University, 
Jena, Germany}
\centerline{\tt kuelshammer@uni-jena.de} %(optional)
\vskip 30pt
%\centerline{\smallit Received: , Accepted: , Published: } % We will fill in the dates
%\vskip 30pt

\baselineskip=10pt
\centerline{\bf Abstract}

%Put your abstract here. Please limit it to half of a page of text.
{\sevenrm
For coprime positive integers $\scriptstyle q$ and $\scriptstyle e$, let $\scriptstyle m(q,e)$ denote the least 
positive integer $\scriptstyle t$ such that there exists a sum of $\scriptstyle t$ powers of $\scriptstyle q$ which 
is divisible by $\scriptstyle e$. We prove that $\scriptstyle m(q,e) \le \lceil e/{\rm ord}_e(q) \rceil$
where $\scriptstyle{\rm ord}_e(q)$ denotes the (multiplicative) order of $\scriptstyle q$ modulo $\scriptstyle e$. We
apply this in order to classify, for any positive integer $\scriptstyle r$, the cases where 
$\scriptstyle m(q,e) \ge {e \over r}$ and $\scriptstyle e > r^4-2r^2$. In particular, we determine all 
pairs $\scriptstyle (q,e)$ such that $\scriptstyle m(q,e) \ge {e \over 6}$. We also investigate in more 
detail the case where $\scriptstyle e$ is a prime power.}

\baselineskip=12.875pt
\parindent=12pt
\vskip 30pt

%Your paper starts here with your first section. See III below for more %information on how to number sections.

\line{\bf 1. Introduction \hfil}
\vskip 10pt
\noindent
Let $q$ and $e$ be coprime positive integers. We denote by $m(q,e)$ the least positive integer $t$ 
such that there exists a sum of $t$ powers of $q$ (repetitions allowed) which is congruent to 0 modulo $e$.
We can, of course, rephrase this in terms of the ring ${\Bbb Z}/e{\Bbb Z}$ of integers modulo $e$, its
group of units $({\Bbb Z}/e{\Bbb Z})^\times$ and the cyclic subgroup $\langle q+e{\Bbb Z} \rangle$ of 
$({\Bbb Z}/e{\Bbb Z})^\times$ generated by the residue class $q+e{\Bbb Z}$. Then $m(q,e)$ is the least
positive integer $t$ such that a sum of $t$ elements in $\langle q+e{\Bbb Z} \rangle$ vanishes in ${\Bbb Z}/e{\Bbb Z}$. 
This interpretation of $m(q,e)$ will be important in the remainder of the paper. 

In order to have an example, note that $4^0 + 4^1 + 4^2 = 21 \equiv 0 \pmod{7}$ but $4^a + 4^b \not\equiv
0 \pmod{7}$ for all nonnegative integers $a,b$; in fact, $4^ a$ and $4^ b$ are each congruent to $1$, 
$4$ or $2$ modulo $7$. 

The function $m(.,.)$ was introduced in [1] in order to bound the Loewy length of certain 
finite-dimensional algebras. Further results on the numbers $m(q,e)$ were established in [2] where 
also tables for small values of $q$ and $e$ can be found. Here we investigate the function $m(.,.)$
in its own right. This will lead to some interesting new results and several open questions. 

Sums of powers also appear, for example, in Waring's problem (cf.\ [6]), in particular in Waring's problem
modulo an integer (cf.\ [8]). However, in Waring's problem one is looking at sums of powers of the form
$$x_1^k + x_2^k + \ldots + x_t^k$$
(i.e.\ the exponent $k$ is fixed) whereas here we are looking at sums of powers of the form
$$q^{a_1} + q^{a_2} + \ldots + q^{a_t}$$
(i.e.\ the base $q$ is fixed). Also, in Waring's problem one wants to write any integer (perhaps modulo $e$)
as a sum of powers whereas we are only interested in writing the single value $0$ modulo $e$ as a sum of powers. 

It is obvious that always $1 \le m(q,e) \le e$. Moreover, we have $m(q,e) = e$ if and only if $q \equiv 
1 \pmod{e}$, and we have $m(q,e) = 1$ if and only if $e = 1$ [1, Example 6.1]. Thus we can ignore these 
trivial cases in the following. We also know that $m(q,e) = 2$ if and only if either $e = 2$, or $n := 
{\rm ord}_e(q)$ is even with $q^{n \over 2} \equiv -1 \pmod{e}$ [1, Lemma 6.3]. Here ${\rm ord}_e(q)$
denotes the (multiplicative) order of $q$ modulo $e$, i.e.\ the least positive integer $s$ such that 
$q^s \equiv 1 \pmod{e}$. Thus ${\rm ord}_e(q)$ is also the order of the cyclic subgroup $\langle 
q+e{\Bbb Z} \rangle$ of $({\Bbb Z}/e{\Bbb Z})^\times$.

When precisely is $m(q,e) = 3$? At present we do not have a good general answer to this question. A
special case is given in [2, Lemma 2.11]. By [1, Lemma 6.2], $m(q,e)$ is always divisible by $e_1 :=
{\rm gcd}(e,q-1)$; this often gives a useful lower bound for $m(q,e)$. In many cases we even have 
$m(q,e) = e_1$; however, we do not have a precise characterization of this phenomenon.

The values $m(q,e)$ are known precisely when $e = 2^k$ for some nonnegative integer $k$ [2, Proposition
2.9]. In fact, in this case we have $m(q,e) = {\rm gcd}(e,q-1) = {e \over {\rm ord}_e(q)}$ if $q \equiv 
1 \pmod{4}$, $m(q,e) = 2$ if $q \equiv -1 \pmod{e}$, and $m(q,e) = 4$ otherwise. In [2, Corollary
2.12], the authors also determined the values $m(q,e)$ when $e$ is a power of a Pierpont prime. Recall 
that a prime $p$ is called a Pierpont prime if $p = 1 + 2^a 3^b$ for suitable nonnegative integers 
$a,b$. Here we will generalize some of the results in [2].

One of our main aims is to improve on the upper bound $m(q,e) \le e$. By [1, Lemma 
6.1], $m(q,e)$ depends only on the residue class of $q$ modulo $e$. Thus, by the remarks above, it is 
enough to concentrate on the case $1 < q < e$. If $e \ge 3$ and $q \equiv -1 \pmod{e}$ then $m(q,e) = 2$
since $q+1$ is divisible by $e$. Thus we can usually even assume $1 < q < e-1$.

Our main result is as follows.
\bigskip\noindent 
{\bf Theorem 1.} {\it Let $q$ and $e$ be coprime positive integers, and set $n := {\rm ord}_e(q)$. Then 
$$m(q,e) \le \left\lceil {e \over n} \right\rceil.$$}
\bigskip\noindent 
For example, Theorem 1 implies that $m(4,7) \le \lceil {7 \over 3} \rceil = 3$. As we observed above, 
we have in fact $m(4,7) = 3$. In particular, it is not always true that $m(q,e) \le {e \over 
{\rm ord}_e(q)}$. It is an open problem to find out when precisely equality holds in the inequality of
Theorem 1. A special case appears in Example 6 below.

Our proof of Theorem 1 will be based on a result of M.~Kneser in additive number theory [4, 6]. We use 
Theorem 1 in order to classify the cases where $m(q,e)$ is ``large''. 
\bigskip\noindent 
{\bf Proposition 2.} {\it Let $q,e$ and $r$ be positive integers such that ${\rm gcd}(q,e) = 1$ and 
$1 < q < e$. We set $e_1 := {\rm gcd}(e,q-1)$.

\noindent (i) If there are positive integers $a,b$ such that 
\medskip 
$(\ast)$ $\quad {\rm gcd}(a,b) = 1$, $b < a \le r$, $ab \le q$ and $e = {a \over b}(q-1)$
\medskip\noindent 
then $m(q,e) = e_1 = {e \over a} \ge {e \over r}$.

\noindent (ii) Conversely, if $m(q,e) \ge {e \over r}$ and $e > r^4 -2r^2$ then there exist positive 
integers $a,b$ satisfying $(\ast)$.}
\bigskip\noindent 
Most likely, the bound $e > r^4-2r^2$ in (ii) is not best possible. We conjecture that it can be 
replaced by the bound $e > 4r^2-4r$. In Section 2 we will use Proposition 2 in order to analyze the 
case $m(q,e) \ge {e \over 6}$ in more detail. This will generalize [2, Proposition 2.5] where the 
case $m(q,e) \ge {e \over 3}$ was considered by a different method. 

In the last section of this paper, we deal with the situation where $e$ is a prime power. Suppose that 
$e = p^k$ where $p$ is an odd prime and $k$ is a positive integer. We reduce the computation of the 
values $m(q,e)$ to the case where $n := {\rm ord}_e(q)$ is not divisible by $p$. Moreover, in this 
situation we show that $m(q,e)$ is the smallest prime divisor of $n$ provided that $k$ is sufficiently 
large. 

In the following, we denote by ${\Bbb Z}$ the ring of integers and by ${\Bbb Q}$ the field of rational 
numbers. Moreover, we denote by ${\Bbb N}$ the set of positive integers and by ${\Bbb P}$ the set of 
prime numbers. We also set ${\Bbb N}_0 := {\Bbb N} \cup \{0\}$. For $a,b \in {\Bbb Z}$, we write $a 
\mid b$ if $a$ divides $b$. For $n \in {\Bbb N}$, we denote the $n$-th cyclotomic polynomial by 
$\Phi_n$. Thus $\Phi_n$ has degree $\varphi(n)$ where $\varphi$ is Euler's totient function. 

The results in this paper have their origin in the first author's bachelor thesis [3] written under 
the guidance of the second author.
\vskip 30pt
\line{\bf 2. General bounds \hfil}
\vskip 10pt
\noindent
We start with the computation of $m(q,e)$ in a special case. We fix coprime positive integers $q$ and $e$
such that $1 < q < e$ and set $e_1 := {\rm gcd}(e,q-1)$ and $n := {\rm ord}_e(q)$. 
\bigskip\noindent 
{\bf Lemma 3.} {\it If $e < e_1^2 + 2e_1$ then $m(q,e) = e_1$. }
\bigskip\noindent 
{\it Proof.} 
%By [1, Lemma 6.2], $m(q,e)$ is always divisible by $e_1$. Thus it suffices to show that 
%$m(q,e) < 2e_1$. 
Our hypothesis implies that $x := {e \over e_1} < e_1 + 2$ and 
$$y := {q-1 \over e_1} \le {e-2 \over e_1} < {e \over e_1} = x \le e_1 + 1.$$
Since ${\rm gcd}(x,y) = 1$ there is an integer $a \in \{0,1,\ldots,y-1\}$ such that $ax \equiv 1 
\pmod{y}$. Then 
$$ax \le (y-1)x \le (y-1)(e_1+1) = ye_1 + y -e_1-1 \le ye_1 - 1 = q-2.$$
Since $e = e_1x = {q-1 \over y}x$ we have 
$$ae = {ax(q-1) \over y} = {ax-1 \over y}q + {q-ax \over y}.$$
If $ax \ne 0$ then ${ax-1 \over y}$ and ${q-ax \over y}$ are nonnegative integers, and thus 
$$m(q,e) \le {ax-1 \over y} + {q-ax \over y} = {q-1 \over y} = e_1 \le m(q,e).$$
If $ax = 0$ then $y=1$ and therefore $e_1 = q-1$. Hence 
$$e < e_1^2+2e_1 = (q-1)^2 + 2(q-1) = q^2-1 = (q-1)q + (q-1).$$
Thus $m(q,e) < 2(q-1) = 2e_1$. By [1, Lemma 6.2], $m(q,e)$ is always divisible by $e_1$. Hence
$m(q,e) = e_1$.  \hfill $\square$
\bigskip\noindent 
We believe that, in fact, the following generalization of Lemma 3 holds:
\bigskip\noindent 
{\bf Conjecture 4.} If $e < (e_1+1)^{k+1}-1$ where $k$ is a positive integer then $m(q,e) \le ke_1$.
\bigskip\noindent 
Explicit computations show that Conjecture 4 is true whenever $1 < q < e < 2050$. The conjecture is also 
true when $e_1 = q-1$. In fact, the $q$-adic expansion 
$$(e_1+1)^{k+1}-1 = q^{k+1}-1 = \sum_{i=0}^k (q-1)q^i$$
shows that $m(q,e) < (k+1)e_1$ if $e < (e_1+1)^{k+1}-1$. Since $m(q,e)$ is divisible by $e_1$ 
this implies: $m(q,e) \le ke_1$. 

For the proof of Theorem 1, we need the following result from additive number theory.
\bigskip\noindent 
{\bf Theorem 5.} (M.~Kneser [7, Theorem 4.3])

\noindent {\it Let $A$ and $B$ be nonempty finite subsets of an abelian group $G$. Then 
$$|A+B| \ge |A+H| + |B+H| - |H|$$
where $H := \{g \in G: g + (A+B) = A+B\}$ is a subgroup of $G$.}
\bigskip\noindent 
Here we use the notation $A+B := \{a+b: a \in A, \, b \in B\}$ and $g+A := \{g+a: a \in A\}$.
The following result is a consequence of Theorem 5.
\bigskip\noindent
{\bf Corollary 6.} {\it Let $R$ be a finite ring, and let $A$ be a subgroup of $R^\times$, its 
group of units. Moreover, let $m(A,R)$ denote the least positive integer $t$ such that there exist 
elements $a_1,\ldots,a_t \in A$ such that $a_1 + \ldots + a_t = 0$. Then 
$$m(A,R) \le \left\lceil {|R| \over |A|} \right\rceil.$$}
\bigskip\noindent 
{\it Proof.} We apply Theorem 5 to the additive group $G$ of $R$ and set $A_1 := A$, $n := |A|$,
$B := A \cup \{0\}$ and $A_m := A_{m-1} + B$ for $1 \ne m \in {\Bbb N}$. Then $A_m$ consists of all 
elements in $R$ which can be written as a sum of $m$ or less elements in $A$. We also set $H_m 
:= \{g \in G: g + A_m = A_m\}$ for $m \in {\Bbb N}$ and claim that the following holds:
\medskip
$(\ast) \quad\quad$ If $0 \notin A_m$ then $|A_m| \ge mn$. 
\medskip\noindent
Certainly, $(\ast)$ holds for $m=1$. Suppose therefore that $m>1$ and that $0 \notin A_m$. Then 
also $0 \notin A_{m-1}$. Arguing by induction, we may assume that $|A_{m-1}| \ge (m-1)n$. 

Assume that $H_m \cap B$ contains a nonzero element $h$. Then $h \in A = A_1 \subseteq A_m$ and 
$h+h \in H_m + A_m = A_m$. Similarly, we obtain $h+h+h \in A_m$. Continuing in this fashion, we
get $kh \in A_m$ for $k \in {\Bbb N}$. In particular, we reach the contradiction $0 = |G|h \in 
A_m$. 

This contradiction shows that, in fact, we have $H_m \cap B = \{0\}$. Then $H_m + B \supseteq
H_m \cup B$ and 
$$|H_m + B| \ge |H_m \cup B| = |H_m| + |B| - |H_m \cap B| 
                     = |H_m| + |B| - 1.$$
Thus Theorem 5 implies:
$$\eqalign{|A_m| &= |A_{m-1} + B| \ge |A_{m-1} + H_m| + |B + H_m| - |H_m| \cr
                 &\ge |A_{m-1}| + |B| - 1 \ge (m-1)n + n = mn.}$$
Thus $(\ast)$ holds. We assume that the corollary is false. Then $m(A,R) > m := \left\lceil {|R| \over 
|A|} \right\rceil$. This means that $0 \notin A_m$. 
Hence $(\ast)$ implies that $|A_m| \ge mn \ge {|R| \over |A|} |A| = |R|$, 
and we have reached the contradiction $0 \in R = A_m$. \hfill $\square$
\bigskip\noindent 
Now Theorem 1 is the special case $R = {\Bbb Z}/e{\Bbb Z}$, $A = \langle q + e{\Bbb Z} \rangle$ of  Corollary 6.                 
\bigskip\noindent 
{\bf Example 7.} Suppose that $q$ is odd and $q>1$, and set $e := 2(q-1)$. Then ${\rm gcd}(q,e) = 1$, 
$n = {\rm ord}_e(q) = 2$ and $e_1 = {\rm gcd}(e,q-1) = q-1$. Since $e < e_1^2+e_1$, Lemma 3 implies 
that $m(q,e) = e_1 = q-1 = {e \over 2} = \lceil {e \over n} \rceil$. Thus, in this case,
the upper bound in Theorem 1 is sharp. 
\bigskip\noindent 
Next we turn to Proposition 2.
\bigskip\noindent 
{\it Proof of Proposition 2.} Let $q,e,r \in {\Bbb N}$ such that ${\rm gcd}(q,e) = 1$ and $1 < q < e$. 

(i) Suppose that there are $a,b \in {\Bbb N}$ such that ${\rm gcd}(a,b) = 1$, $b < a \le r$, $ab \le q$
and $e = {a \over b}(q-1)$. Since $b$ divides $a(q-1)$ and ${\rm gcd}(a,b) = 1$ we conclude that $b$ 
divides $q-1$. Thus 
$$e_1 = {\rm gcd}(e,q-1) = {q-1 \over b}{\rm gcd}(a,b) = {q-1 \over b}.$$
Moreover, $ab \le q$ implies that $a \le {q \over b} = {q-1 \over b} + {1 \over b} \le e_1+1$. Hence
$$e = a {q-1 \over b} \le (e_1+1)e_1 = e_1^2+e_1 < e_1^2+2e_1.$$
Thus Lemma 3 implies that $m(q,e) = e_1 = {q-1 \over b} = {e \over a} \ge {e \over r}$ since $a \le r$.

(ii) Suppose that $m(q,e) \ge {e \over r}$ and $e > r^4-2r^2$. If $r=1$ then $m(q,e) = e$. Hence $q 
\equiv 1 \pmod{e}$ which is impossible. Thus we must have $r \ge 2$.

If $e < e_1^2+2e_1$ then Lemma 3 implies that $m(q,e) = e_1$. We set $a := {e \over e_1}$ and $b := 
{q-1 \over e_1}$. Then ${\rm gcd}(a,b) = 1$ and $e = ae_1 = {a \over b}(q-1)$. Since $q < e$ we have $b < {e-1 \over e_1} < {e \over e_1} = a$.
Also, $e_1 = m(q,e) \ge {e \over r} = {ae_1 \over r}$ implies that $a \le r$. 
Assuming $q < ab$ we get the contradiction 
$$r^4-2r^2 < e = {a \over b} (q-1) < {a \over b} (ab-1) = a(a-{1 \over b}) < a^2 \le r^2.$$
It remains to consider the case $e \ge e_1^2+2e_1$. Theorem 1 implies that ${e \over r} \le m(q,e) \le 
\lceil {e \over n} \rceil$. If $n>r$ then ${e \over r} \le \lceil {e \over 
r+1} \rceil < {e \over r+1} + 1 = {e+r+1 \over r+1}$, i.e.\ $er+e < er+r^2+r$. This leads to the contradiction
$r^4-2r^2 < e < r^2+r$.

Thus we must have $n \le r$. Hence [1, Lemma 6.2] implies that $m(q,e) \le re_1$. Therefore
$e_1^2+2e_1 \le e \le m(q,e)r \le r^2e_1$, so that $e_1 \le r^2-2$. Now we have reached the 
contradiction $e \le r^2e_1 \le r^4-2r^2$. \hfill $\square$
\bigskip\noindent 
Next we analyze the situation for $m(q,e) \ge {e \over 6}$ in more detail.
\bigskip\noindent
{\bf Corollary 8.} {\it Suppose that $1<q<e-1$.  Then 
$$m := m(q,e) \ge {e \over 6}$$
if and only if one of the following holds:

\noindent (i) $q \ge a+1$, ${\rm gcd}(a,q) = 1$ and $e = a(q-1)$ for some $a \in \{2,3,4,5,6\}$ 

\noindent (in which case $m = {e \over a} = q-1$);
%\noindent (i) $q \ge 3$, ${\rm gcd}(2,q) = 1$ and $e = 2(q-1)$ \hfil (in which case $m = {e  \over 2} =
%q-1$);
%\noindent (ii) $q \ge 4$, ${\rm gcd}(3,q) = 1$ and $e = 3(q-1)$ \hfil (in which case $m = {e \over 3} =
%q-1$);

\noindent (ii) $q \ge 7$, ${\rm gcd}(6,q) = 1$ and $e = {3 \over 2}(q-1)$ \hfill (in which case 
$m = {e \over 3} = {q-1 \over 2}$);
%\noindent (iv) $q \ge 5$, ${\rm gcd}(2,q) = 1$ and $e = 4(q-1)$ \hfil (in which case $m = {e \over 4} = 
%q-1$);

\noindent (iii) $q \ge 13$, $q \equiv 1 \pmod{6}$ and $e = {4 \over 3}(q-1)$ \hfill 

\noindent (in which case $m = {e \over 4} = {q-1 \over 3}$);
%\noindent (vi) $q \ge 6$, ${\rm gcd}(5,q) = 1$ and $e = 5(q-1)$ \hfil (in which case $m = {e \over 5} =
%q-1$);

\noindent (iv) $q \ge 7$, ${\rm gcd}(10,q) = 1$ and $e = {5 \over 2}(q-1)$ \hfill (in which case $m =
{e \over 5} = {q-1 \over 2}$);

\noindent (v) $q \ge 7$, ${\rm gcd}(5,q) = 1$, $q \equiv 1 \pmod{3}$, $e = {5 \over 3}(q-1)$ 

\noindent (in which case $m = {e \over 5} = {q-1 \over 3}$);

\noindent (vi) $q \ge 13$, ${\rm gcd}(5,q) = 1$, $q \equiv 1 \pmod{4}$, $e = {5 \over 4}(q-1)$  

\noindent (in which case $m = {e \over 5} = {q-1 \over 4}$);
%\noindent (x) $q \ge 7$, ${\rm gcd}(6,q) = 1$ and $e = 6(q-1)$ \hfil (in which case $m = {e \over 6} =
%q-1$);

\noindent (vii) $q \ge 16$, ${\rm gcd}(6,q) = 1$, $q \equiv 1 \pmod{5}$, $e = {6 \over 5}(q-1)$ 

\noindent (in which case $m = {e \over 6} = {q-1 \over 5}$);

\noindent (viii) (cases where $m(q,e) = 2$; cf.\ [1, Lemma 6.3])

\noindent $e = 5$, $q \in \{2,3\}$; $e = 7$, $q \in \{3,5\}$; $e = 9$, $q \in \{2,5\}$; $e = 10$, $q  \in \{3,7\}$; $e = 11$, $q \in \{2,6,7,8\}$;

\noindent (ix) (cases where $n = 2$ and $m(q,e) = 2e_1 > 2$; cf.\ [2, Proposition 2.7])

\noindent $(e,q) \in \{(8,3),\,(15,4),\,(16,7),\,(21,13),\,(24,5),\,(24,11),\,(33,10),\,(35,6),$

\noindent $(40,29),(45,26),(48,7),(55,21),(63,8),(77,43),(80,9),(99,10),(120,11)\}$;

\noindent (x) $e = 7$, $q \in \{2,4\}$, $m=3$; $e=11$, $q \in \{3,4,5,9\}$, $m=3$; $e = 13$, $q \in 
\{3,9\}$, $m=3$; $e=14$, $q \in \{9,11\}$, $m=4$; $e=15$, $q \in \{2,8\}$, $m=4$; $e = 16$, $q \in 
\{3,11\}$, $m=4$; $e = 20$, $q \in \{3,7\}$, $m=4$; $e = 22$, $q \in \{3,5,9,15\}$, $m=4$; $e=26$,
$q \in \{3,9\}$, $m=6$; $e = 48$, $q \in \{5,29\}$, $m=8$.}
\bigskip\noindent 
{\it Proof.} Suppose first that $e \le 6^4-2 \cdot 6^2 = 1224$. Then a direct computation of $m(q,e)$
for all relevant $q$ and $e$ shows that the result holds in this case. Thus, in the following, we may
assume that $e > 1224$. 

If $m(q,e) \ge {e \over 6}$ then, by Proposition 2, there exist coprime integers $a,b \in 
{\Bbb N}$ such that $b < a \le 6$, $ab \le q$ and $e = {a \over b}(q-1)$. This leads to one of the 
cases (i), $\ldots$, (vii).

Now suppose, conversely, that we are in one of the cases (i), $\ldots$, (vii). Then $e = {a \over b}
(q-1)$ with coprime $a,b \in {\Bbb N}$. Moreover, we have $b < a \le 6$ and $ab \le q$ in all cases. 
Thus, again by Proposition 2, we have $m(q,e) = {\rm gcd}(e,q-1) = {e \over a} = {q-1 \over b} \ge 
{e \over 6}$. \hfill $\square$
\bigskip\noindent 
Let us consider the cases (viii), (ix) and (x) as exceptional. Then, for $c=2,3,4,5,6$, the ``largest''
exceptional case for the inequality $m(q,e) \ge {e \over c}$ occurs when $e = 4c^2-4c$ and $q = 2c-1$.

%In a certain sense, the cases (viii) and (ix) are well understood. The ``largest'' case in 
%(ix) is much smaller than the ``largest'' case in (x).
\vskip 30pt
\line{\bf 3. When $e$ is a prime power \hfil}
\vskip 10pt
\noindent
In this section, we investigate the numbers $m(q,e)$ in the case where $e$ is a power of a prime  $p$. 
Suppose that $p \nmid q \in {\Bbb N}$. Then [1, Lemma 6.1] implies that
$$m(q,p) \le m(q,p^2) \le m(q,p^3) \le \ldots$$
and 
$$m(q,p^{j+k}) \le m(q,p^j)m(q,p^k)$$
for $j,k \in {\Bbb N}$. Moreover, Theorem 1 implies that 
$$m(q,p^k) \le \left\lceil {p^k \over {\rm ord}_{p^k}(q)} \right\rceil$$
for $k \in {\Bbb N}$. 
It is easy to see that the sequence $(p^k/ {\rm ord}_{p^k}(q))_{k=1}^\infty$ is bounded if $q \ne 1$
(cf.\ the proof of Proposition 9 below). Thus the sequence $(m(q,p^k))_{k=1}^\infty$ is also 
bounded. In the following, we will investigate this sequence and its limit in more detail.
For example, [2, Proposition 2.17] implies that the sequence $(m(9,11^k))_{k=1}^\infty$ is given by
$(3,5,5,\ldots)$ whereas the sequence $(11^k / {\rm ord}_{11^k}(9))_{k=1}^\infty$ is given by
$({11 \over 5}, {121 \over 5}, {121 \over 5}, \ldots)$.

As we mentioned in the introduction, we know all values of the form $m(q,2^k)$.
Thus, in the following, we will usually assume that $p>2$. Then $({\Bbb Z}/p^k{\Bbb Z})^\times$ is a 
cyclic group of order $\varphi(p^k) = p^{k-1}(p-1)$. Moreover, by [1, Lemma 6.1], the value $m(q,p^k)$ 
depends only on $p^k$ and ${\rm ord}_{p^k}(q)$, rather than on $p^k$ and $q$. Our first result in this section 
is an application of Hensel's Lemma; it generalizes [2, Lemma 2.14]. 
\bigskip\noindent 
{\bf Proposition 9.} {\it Let $2 \ne p \in {\Bbb P}$, let $k,q \in {\Bbb N}$ such that $p \nmid q$, and write 
${\rm ord}_{p^k}(q) = p^i d$ where $i,d \in {\Bbb N}_0$ and $d \mid p-1$. 
If $i > 0$ then ${\rm ord}_{p^{k-1}}(q) = p^{i-1}d$ and $m(q,p^k) = m(q,p^{k-1})$. }
\bigskip\noindent 
{\it Proof.} Suppose that $i>0$, so that $k \ge i+1 \ge 2$. The kernel $K$ of the canonical 
epimorphism 
$$f: ({\Bbb Z}/p^k{\Bbb Z})^\times \longrightarrow ({\Bbb Z}/p^{k-1}{\Bbb Z})^\times, \quad 
x + p^k{\Bbb Z} \longmapsto x+p^{k-1}{\Bbb Z},$$
has order $p$. Since $({\Bbb Z}/p^k{\Bbb Z})^\times$ is cyclic, $K$ is contained in $\langle 
q+p^k{\Bbb Z} \rangle$. Thus $\langle q+p^{k-1}{\Bbb Z} \rangle = f(\langle q+p^k{\Bbb Z} \rangle)$ has
order $p^{i-1}d$, i.e.\ ${\rm ord}_{p^{k-1}}(q) = p^{i-1}d$. 

We set $c := {p-1 \over d}$. Then every element in $\langle q+p^{k-1}{\Bbb Z} \rangle$ is a $p^{k-1-i}c
$-th power in $({\Bbb Z}/p^{k-1}{\Bbb Z})^\times$. Thus there is a relation of the form $x_1^c + 
\ldots + x_t^c \equiv 0 \pmod{p^{k-1}}$ where $t = m(q,p^{k-1})$ and each $x_j + p^{k-1}{\Bbb Z}$ is 
a $p^{k-1-i}$-th power in $({\Bbb Z}/p^{k-1}{\Bbb Z})^\times$. Since $c \not\equiv 0 \pmod{p}$, 
Hensel's Lemma implies that we also have a relation of the form $y_1^c + \ldots + y_t^c \equiv 0 
\pmod{p^k}$ where $y_1,\ldots,y_t \in {\Bbb Z}$ satisfy $y_j \equiv x_j \pmod{p^{k-1}}$ for $j=1,\ldots,
t$. Thus 
$$y_j^{p^{i-1}(p-1)} \equiv x_j^{p^{i-1}(p-1)} \equiv 1 \pmod{p^{k-1}}$$
and 
$$1 \equiv y_j^{p^i (p-1)} \equiv (y_j^c)^{dp^i} \pmod{p^k},$$
i.e.\ $y_j^c \in \langle q+p^k{\Bbb Z} \rangle$ for $j=1,\ldots,t$. This shows that 
$$m(q,p^k) \le t = m(q,p^{k-1}) \le m(q,p^k),$$
and the result follows. \hfill $\square$
\bigskip\noindent 
Proposition 9 reduces the computation of the values $m(q,p^k)$ to the situation where ${\rm 
ord}_{p^k}(q)$ divides $p-1$. Note that in this case we have ${\rm ord}_{p^k}(q) = {\rm 
ord}_p(q)$.

We start again from a somewhat different angle. Let $2 \ne p \in {\Bbb P}$, let $q \in {\Bbb N}$ such that 
$p \nmid q$, and suppose that $n := {\rm ord}_p(q)$ divides $p-1$. The case $n=1$ is trivial: we have 
$q \equiv 1 \pmod{p}$ and $m(q,p^k) = {\rm gcd}(p^k,q-1) = p^k/{\rm ord}_{p^k}(q)$ for $k \in {\Bbb N}$, by [2,
Proposition 2.10]. Thus, in the following, we may assume that $n>1$. 

We denote by $p^w$ the maximal power of $p$ dividing $q^n-1$. Then $w > 0$ since $p \mid q^n-1$. (We also have
$w = v_p(q^n-1)$ where $v_p$ denotes the $p$-adic valuation of ${\Bbb Q}$.) Thus $q^n \equiv 1 \pmod{p^k}$ for 
$k=1,\ldots,w$, but $q^n \not\equiv 1 \pmod{p^{k+1}}$. This means that
$$n = {\rm ord}_p(q) = {\rm ord}_{p^2}(q) = \ldots = {\rm ord}_{p^w}(q) < {\rm ord}_{p^{w+1}}(q) < \ldots $$
where ${\rm ord}_{p^{w+j}}(q) = np^j$ for $j \in {\Bbb N}$. Furthermore, Proposition 9 implies that
$$m(q,p) \le m(q,p^2) \le \ldots \le m(q,p^w) = m(q,p^{w+1}) = \ldots .$$
Thus $\lim_{k \to \infty} m(q,p^k) = m(q,p^ w) \le r$ where $r$ denotes the smallest prime divisor of $n$, by
[1, Lemma 6.4]. Hence $r \mid n \mid p-1$. 

At this point we recall that (cf.\ [1, Remark 6.1]) $m(q,p^k) = 2$ if and only if ${\rm ord}_{p^k}(q)$ is even.
Moreover, [2, Lemma 2.11] shows that $m(q,p^k) = 3$ if $3$ is the smallest prime divisor
of ${\rm ord}_{p^k}(q)$ and $p > 3$. Hence we will often assume that $r \ge 5$. 
\bigskip\noindent 
{\bf Proposition 10.} {\it Let $2 \ne p \in {\Bbb P}$, let $1 \ne n \in {\Bbb N}$ such that $n \mid p-1$, 
and let $r$ be the smallest prime divisor of $n$. 
Then there exist $K,Q \in {\Bbb N}$ such that $p \nmid Q$, ${\rm ord}_{p^K}(Q) = n$ and $m(Q,p^K) = r$.}
\bigskip\noindent 
{\it Proof.} We denote by ${\Bbb Z}_p$ the ring of $p$-adic integers. It is well-known that ${\Bbb Z}_p$
contains a primitive $n$-th root of unity $\zeta$. The main result of [5] shows that $\zeta^ {i_1} +
\ldots + \zeta^{i_s} \ne 0$ in ${\Bbb Z}_p$ for $s=1,\ldots,r-1$ and $i_1,\ldots,i_s \in \{1,\ldots,n\}$.
(The authors of [5] work with roots of unity in ${\Bbb C}$, but it is clear that their results also
hold in ${\Bbb Z}_p$.) Thus, for any such numbers $s,i_1,\ldots,i_s$ there exists a positive integer 
$k = k(s;i_1,\ldots,i_s)$ such that $\zeta^{i_1} + \ldots + \zeta^{i_s} \not\equiv 0 \pmod{p^k{\Bbb Z}_p}$.
We set 
$$K := \max \{k(s;i_1,\ldots,i_s): 1 \le s < r, \; 1 \le i_j \le n \hbox{ for } j=1,\ldots,s \}.$$
Then $\zeta^{i_1} + \ldots + \zeta^{i_s} \not\equiv 0 \pmod{p^K{\Bbb Z}_p}$ for $s=1,\ldots,r-1$ and 
$i_1,\ldots,i_s \in \{1,\ldots,n\}$. Now we use the ring isomorphism
$${\Bbb Z}/p^K{\Bbb Z} \longrightarrow {\Bbb Z}_p / p^K{\Bbb Z}_p, \quad x+p^K{\Bbb Z} \longmapsto
x+p^K{\Bbb Z}_p.$$
Let $Q \in {\Bbb N}$ such that $Q + p^K{\Bbb Z}$ maps to $\zeta + p^K{\Bbb Z}_p$ under this isomorphism.
Since $\zeta + p^K{\Bbb Z}_p$ is an element of order $n$ in $({\Bbb Z}_p/p^K{\Bbb Z}_p)^\times$, 
$Q + p^K{\Bbb Z}$ is an element of order $n$ in $({\Bbb Z}/p^K{\Bbb Z})^\times$. This means that $p \nmid Q$
and ${\rm ord}_{p^K}(Q) = n$. Moreover, we have $Q^{i_1} + \ldots + Q^{i_s} \not\equiv 0 \pmod{p^K
{\Bbb Z}}$ for $s=1,\ldots,r-1$ and $i_1,\ldots,i_s \in \{1,\ldots,n\}$. This implies that $m(Q,p^K)
\ge r$. Since also $m(Q,p^K) \le r$ the result follows. \hfill $\square$
\bigskip\noindent 
Let $2 \ne p \in {\Bbb P}$, let $n \in {\Bbb N}$ such that $1 \ne n \mid p-1$, and let $r$ be the 
smallest prime divisor of $n$. Then, by Proposition 10, there is $k \in {\Bbb N}$ such that 
$$m(q,p) \le m(q,p^2) \le \ldots \le m(q,p^k) = r,$$
for suitable $q \in {\Bbb N}$ such that $p \nmid q$ and ${\rm ord}_{p^k}(q) = n$. Also, [1, Lemma 6.1]
implies that 
$$m(q,p^j) \le m(q^i,p^j)$$
for $i,j \in {\Bbb N}$. In 
addition, we recall that $m(q,p^j) = 2$ whenever $r = 2$, and $m(q,p^j) = 3$ whenever $r = 3$. 

The results above show that, for a fixed odd prime $p$, the knowledge of all values $m(q,p^k)$ requires
only the computation of finitely many values. Below we will give several examples.
It would be nice to have an upper bound for the number $K$ in Proposition 10.

Our next result generalizes [2, Proposition 2.5]. Here the condition $q^n-1 \equiv 0 \pmod{e}$ coming
from the hypothesis ${\rm ord}_e(q) = n$ is replaced by the stronger condition $\Phi_n(q) \equiv 0
\pmod{e}$. We will see later that, at least in favourable situations (cf.\ Propositions 14 and 15), 
the values $m(q,p^k)$ depend more on ${\rm ord}_{p^k}(q)$ than on $p$, $k$ and $q$.
\bigskip\noindent 
{\bf Proposition 11.} {\it Let $1 \ne n \in {\Bbb N}$. Then there are only finitely many $e \in {\Bbb N}$
such that, for some $q \in {\Bbb N}$, we have $e \mid \Phi_n(q)$ and $m(q,e) < {n \over n - \varphi(n)}$.}
\bigskip\noindent 
The proof will show that - at least in principle - the finitely many possibilities for $e$ can be
computed. In [2], this has been done for $n = 5$ and $n = 7$. 
\bigskip\noindent
{\it Proof.} Let $q,e \in {\Bbb N}$. If $e$ divides $\Phi_n(q)$ then $e$ also divides $q^n-1$; in 
particular, we have ${\rm gcd}(q,e) = 1$. Thus $m := m(q,e)$ is defined. 

Now suppose, in addition, that $m < {n \over n-\varphi(n)}$. Then there are $i_1,\ldots,i_m \in {\Bbb 
N}_0$ such that 
$$i_1 \le \ldots \le i_m \quad {\rm and} \quad q^{i_1} + \ldots + q^{i_m} \equiv 0 \pmod{e}.$$
Here we may assume that $i_m$ is as small as possible. Then $i_m < n$, and $i_1 = 0$ since otherwise we 
can cancel powers of $q$. We claim that $i_m < \varphi(n)$.

If $i_{l+1} - i_l \le n-\varphi(n)$ for $l=1,\ldots,m-1$ then 
$$\eqalign{i_m &= (i_m -i_{m-1}) + \ldots + (i_2-i_1) \le (m-1)(n-\varphi(n)) \cr
&= m (n-\varphi(n))-n+\varphi(n) < n-n+\varphi(n) = \varphi(n).\cr}$$
If $i_{l+1}-i_l > n-\varphi(n)$ for some $l \in \{1,\ldots,m-1\}$ then we multiply our relation by 
$q^{n-i_{l+1}}$ and reorder the summands to obtain a new relation
$q^{j_1} + \ldots + q^{j_m} \equiv 0 \pmod{e}$ where $0 \le j_1 \le \ldots \le j_m = 
i_l + n - i_{l+1} < \varphi(n).$
Thus our claim is proved. 

The polynomial $g := X^{i_1} + \ldots + X^{i_m} \in {\Bbb Q}[X]$ is not divisible
by $\Phi_n$ since ${\rm deg}(\Phi_n) = \varphi(n) > i_m = {\rm deg}(g)$. Since $\Phi_n$ is irreducible in 
${\Bbb Q}[X]$ there are $a,b \in {\Bbb Q}[X]$ such that $ag+b\Phi_n = 1$. We choose $d \in {\Bbb N}$ minimal
such that $da,db \in {\Bbb Z}[X]$. Then 
$$d = (da)(q) \cdot g(q) + (db)(q) \cdot \Phi_n(q) \equiv 0 \pmod{e}, \quad {\rm i.e.} \quad e \mid d.$$
Since there are only finitely many choices for $i_1,\ldots,i_m$ the result follows. \hfill $\square$
\bigskip\noindent 
{\bf Remark 12.} (i) Note that, for $1 \ne n \in {\Bbb N}$, we have 
$${n \over n-\varphi(n)} = {n' \over n'-\varphi(n')} \quad {\rm where} \quad n' := {\rm rad}(n).$$
(Recall that ${\rm rad}(n) = r_1 \ldots r_s$ whenever the prime factorization of $n$ is $n = r_1^{a_1}
\ldots r_s^{a_s}$.) Thus, in the following, we can often assume that $n$ is squarefree.

\noindent (ii) It is easy to see that always 
$${n \over n-\varphi(n)} \le r$$
where $r$ is the smallest prime divisor of $n$. 

\noindent (iii) It is also easy to see that ``many'' integers $n$ (in particular, all prime powers) satisfy 
$${n \over n - \varphi(n)} > r-1.$$
%Suppose that $n = r_1 \ldots r_s$ where $r_1,\ldots,r_s \in {\Bbb P}$ and $r := r_1 < 
%\ldots < r_s$. We set $R := r_2 \ldots r_s$. A straightforward computation shows that 
%$$ {n \over n-\varphi(n)} > r-1 \hbox{ if and only if } {\varphi(R) \over R} > {r^2 - 2r \over 
%r^2 -2r + 1}.$$
%Moreover, this holds if and only if $r$ is arbitrary, $r_2 > (r-1)^2$, $r_3 > {(r_2-1)(r-1)^2 \over 
%r_2 - (r-1)^2}$, etc. Thus ``most'' integers $n$ (and all prime powers $n$) satisfy ${n \over n - 
%\varphi(n)} > r-1$. 
\bigskip\noindent 
Proposition 11 has the following consequence.
\bigskip\noindent 
{\bf Corollary 13.} {\it Let $1 \ne n \in {\Bbb N}$. Then there exist only finitely many pairs 
$(p,k) \in {\Bbb P} \times {\Bbb N}$ such that, for some $q \in {\Bbb N} \setminus p{\Bbb N}$,
we have $n \mid p-1$, ${\rm ord}_{p^k}(q) = n$ and $m(q,p^k) < {n \over n - \varphi(n)}$.}
%the following set is finite:
%$$\{(p,k) \in {\Bbb P} \times {\Bbb N}\colon \exists q \in {\Bbb N}\setminus p{\Bbb N}\colon 
%n \mid p-1,  
%{\rm ord}_{p^k}(q)= n, m(q,p^k) < {n \over n-\varphi(n)}\}.$$}
\bigskip\noindent 
{\it Proof.} Let $p \in {\Bbb P}$ such that $n \mid p-1$, and let $q,k \in {\Bbb N}$ such that $p 
\nmid q$, ${\rm ord}_{p^k}(q) = n$ and $m(q,p^k) < {n \over n-\varphi(n)}$. Then
$$e := p^k \mid q^n-1 = \prod_{d \mid n} \Phi_d(q).$$
If $p \mid \Phi_d(q)$ for some $d \mid n$ then also $p \mid q^d-1$, i.e.\ $q^d \equiv 1 \pmod{p}$. 
Thus ${\rm ord}_e(q^d)$ is a power of $p$. On the other hand, we have ${\rm ord}_e(q^d) \mid {\rm 
ord}_e(q) = n \mid p-1$. We conclude that ${\rm ord}_e(q^d) = 1$, so that $q^d \equiv 1 \pmod{e}$.
This implies: $d = n$. 

Now, since $p \nmid \Phi_d(q)$ for every proper divisor $d$ of $n$, we must have $e \mid \Phi_n(q)$. 
Thus $e = p^k$ is contained in the finite set described in Proposition 11. \hfill $\square$
\bigskip\noindent 
Corollary 13 implies that, in many cases, the number $K$ in Proposition 10 will be $1$. We illustrate 
this fact by the next two propositions.
\bigskip\noindent 
{\bf Proposition 14.} {\it Let $p \in {\Bbb P}$ such that $5 \mid p-1$, and let $k,q \in {\Bbb N}$ such 
that $p \nmid q$ and ${\rm ord}_{p^k}(q) = 5$. Then 
$$m(q,p^k) = \cases{3, &if $p=11$ and $k=1$, \cr
4, &if $p=61$ and $k=1$, \cr
5, &otherwise.\cr}$$}
\bigskip\noindent 
{\it Proof.} Let $n := 5$. Then [2, Example 2.26] shows that the finite set of Corollary 13 is 
$\{(11,1),(61,1)\}$. If $(p,k)$ is not in this set, then $m(q,p^k) \ge {5 \over 5-4} = 5$, by 
Corollary 13. Thus $m(q,p^k) = 5$ in these cases. 

For $(p,k) \in \{(11,1),(61,1)\}$ the value $m(q,p^k)$ is obtained by a direct computation. \hfill $\square$
\bigskip\noindent 
In a similar way, we can handle the case $n := {\rm ord}_{p^k}(q) = 7$.
\bigskip\noindent 
{\bf Proposition 15.} {\it Let $p \in {\Bbb P}$ such that $7 \mid p-1$, and let $k,q \in {\Bbb N}$
such that $p \nmid q$ and ${\rm ord}_{p^k}(q) = 7$. Then 
$$m(q,p^k) = \cases{3, &if $p=43$ and $k=1$,\cr
                    4, &if $p \in \{29,71,547\}$ and $k=1$,\cr
                    5, &if $p \in \{113,197, 421, 463\}$ and $k=1$,\cr
                    6, &if $p \in \{211,379,449,757,2689\}$ and $k=1$,\cr
                    7, &otherwise.\cr}$$}
\bigskip\noindent   
Next we consider some cases where $n := {\rm ord}_{p^k}(q)$ is a larger prime.
\bigskip\noindent 
{\bf Example 16.} (i) Let $n := 11$. The following table contains the sequence $(m(q,p),m(q,p^2), \ldots)$ for some primes $p$ where $p \nmid q$ and ${\rm ord}_{p^j}(q) = 11$ for $j=1,2,\ldots$:
\medskip\noindent 
$$\matrix{\hfill p = 23:            & (3,5,9,9,11)\hfill & \quad & \hfill p = 67:  & (4,8,11)\hfill \cr
          \hfill p = 89:            & (4,9,11)\hfill     & \quad & \hfill p = 199: & (6,11)\hfill   \cr
          \hfill p \in \{353,397\}: & (5,11)\hfill       & \quad &                 &                 \cr}$$
\medskip\noindent 
(ii) Replacing $n=11$ by $n=13$, we have the following examples:
\medskip\noindent 
$$\matrix{\hfill p=53:  & (3,7,12,13)\hfill & \quad & \hfill p \in \{79,157\}:      & (4,8,12,13)\hfill \cr
          \hfill p=131: & (4,8,13)\hfill    & \quad & \hfill p=313:                 & (5,10,13)\hfill   \cr
          \hfill p=521: & (7,13)\hfill      & \quad & \hfill p \in \{547,677,937\}: & (5,13)\hfill      \cr
          \hfill p=911: & (6,13)\hfill      & \quad &                               &                   \cr}$$ 
\medskip\noindent
(iii) For $n=17$, we have the following examples:
\medskip\noindent 
$$\matrix{\hfill p=239: & (3,9,15,17)\hfill  & \quad & \hfill p=307: & (4,9,14,17)\hfill \cr
          \hfill p=409: & (5,10,15,17)\hfill & \quad & \hfill p=613: & (5,10,17)\hfill   \cr
          \hfill p=919: & (5,12,17)\hfill    & \quad & \hfill p=953: & (4,11,17)\hfill   \cr}$$
\medskip\noindent 
(iv) For $n=19$, we obtain:
\medskip\noindent 
$$\matrix{\hfill p=229: & (5,8,11,19)\hfill & \quad & \hfill p=571: & (4,9,16,19)\hfill \cr
          \hfill p=761: & (5,7,17,19)\hfill &       &               &                   \cr}$$
\bigskip\noindent 
Finally, we present some examples where $n := {\rm ord}_{p^k}(q)$ is composite.
\bigskip\noindent 
{\bf Example 17.} (i) (Cases where the smallest prime divisor of $n$ is $r=5$)
\medskip\noindent 
$(m(q,p),m(q,p^2)) = (3,5)$ for $(p,n) \in \{(71,35), (101,25), (131,65),$

\noindent $(211,35),(281,35), (521,65),
(571,95), (631,35), (911,35)\};$
\medskip 
\noindent $(m(q,p),m(q,p^2)) = (4,5)$ for $(p,n) \in \{(421,35), (491,35), (701,35),$ 

\noindent $(761,95),(911,65), (1051,35), (1471,35), (2311,35), (2521,35), (2591,35),$ 

\noindent $(2731,35), (3221,35), (3361,35),(3571,35), (3851,35)\};$
\medskip
\noindent $m(q,p) = 5$ for $(p,n) \in \{(1151,25), (1201,25), (1301,25), (1801,25),$ 

\noindent $(2381,35), (2801,35), (2861,55), (3011,35)\}.$
\medskip\noindent 
(ii) (Cases where the smallest prime divisor of $n$ is $r=7$)
\medskip\noindent 
%$p=71, \; n = 35: (3,5);$ \hfill $p=101, \; n = 25: (3,5);$
%\noindent $p = 131, \; n = 65: (3,5);$ \hfill 
$$\matrix{\hfill p = 239 & \hfill n = 119: & (3,4,6,7)\hfill & \quad  & p = 547 & n = 91: & (3,4,7) \cr
          \hfill p = 911 & \hfill n = 91:  & (4,6,7)\hfill   &        &         &         &         \cr}$$ 
%\noindent $p = 571, \; n = 95: (3,5);$ \hfill $p = 761, \; n = 95: (4,5);$
%\noindent $p = 911, \; n = 35: (3,5);$ \hfill $p = 911, \; n = 65: (4,5);$

%\hfill $p = 911, \; n = 455: (3,4,5).$
\vskip 30pt
\line{\bf References \hfil}
\vskip 10pt
\noindent 
\item{[1]} T.~Breuer, L.~H\'ethelyi, E.~Horv\'ath and B.~K\"ulshammer, The Loewy structure of certain
fixpoint algebras, Part I, {\it J.\ Algebra} {\bf 558} (2020), 199-220
\smallskip\noindent 
\item{[2]} T.~Breuer, L.~H\'ethelyi, E.~Horv\'ath and B.~K\"ulshammer, The Loewy structure of certain
fixpoint algebras, Part II, {\it Intern.\ Electr.\ J.\ Algebra} {\bf 30} (2021), 16-65
\smallskip\noindent 
\item{[3]} L.~Jacob, \"Uber Summen von Potenzen einer nat\"urlichen Zahl, Bachelor thesis, Jena 2020
\smallskip\noindent 
\item{[4]} M.~Kneser, Absch\"atzungen der asymptotischen Dichte von Summenmengen, {\it Math.\ Z.} 
{\bf 58} (1953), 459-484
\smallskip\noindent 
\item{[5]} T.~Y.~Lam and K.~H.~Leung, On vanishing sums of roots of unity, {\it J.\ Algebra} {\bf 224}
(2000), 91-109
\smallskip\noindent 
\item{[6]} M.~B.~Nathanson, Additive number theory: The classical bases, Springer-Verlag, New York 
1996
\smallskip\noindent
\item{[7]} M.~B.~Nathanson, Additive number theory: Inverse problems and the geometry of sumsets, 
Springer-Verlag, New York 1996
\smallskip\noindent
\item{[8]} C.~Small, Solution of Waring's problem mod $n$, {\it Amer.\ Math.\ Monthly} {\bf 84} (1977),
356-359
\end

%III. In the body of the paper, please follow A, B, and C below.

%A. Label and number your (sub)sections as follows.

%i). "\line{\bf N. Section Name \hfil}", where N=1,2,... corresponds to the section number.

%and/or

%ii). "\line{\bf N.x Subsection Name \hfil}", where N=1,2,3,... and x=1,2,3,... correspond to the subsection number.

%Follow each of the above "line" commands by:
%\vskip 10pt
%\noindent

%B. After the end of each section and before the beginning of the next section please put the line:

%\vskip 30pt

C. Please label your Definitions, Lemmas, Theorems, etc, in the format "{\bf Theorem x}", where x is how you see fit. Further, denote the start of your proofs by "{\it Proof.}".

IV. The references and acknowledgments sections of the paper should not be numbered sections. Hence, you should use "\line{\bf References \hfil}" and/or "\line{\bf Acknowledgments \hfil}." Please label your references numerically. Furthermore, please use the AMS abbreviations for journals (see this link).

\end  

\bigskip\noindent 
{\it Proof of Theorem 1.} We apply Theorem 5 to the abelian group $G := {\Bbb Z}/e{\Bbb Z}$. The 
multiplicative group $A_1 := \langle q+e{\Bbb Z} \rangle$ generated by $q+e{\Bbb Z}$ has order $n := 
{\rm ord}_e(q)$. We set $B := A_1 \cup \{0\}$ and $A_m := A_{m-1}+B$ for $m>1$. Then $A_m$ consists 
of all elements in $G$ which can be written as sums of $m$ or less powers of $q+e{\Bbb Z}$. We also
define $H_m := \{g \in G: g+A_m = A_m\}$ for $m \in {\Bbb N}$ and claim that the following is true:
\medskip
$(\ast)$ If $0+e{\Bbb Z} \notin A_m$ then $|A_m| \ge mn$.
\medskip\noindent 
Certainly, $(\ast)$ holds for $m=1$. Suppose that $m>1$. If $0+e{\Bbb Z} \notin A_m$ then also $0+
e{\Bbb Z} \notin A_{m-1}$. Arguing by induction, we may assume that $|A_{m-1}| \ge (m-1)n$. 

Assume that $H_m \cap B$ contains a nonzero element $h$. Then $h \in A_1 \subseteq A_m$ and $h+h \in 
H_m + A_m = A_m$. Similarly, we obtain $h+h+h \in A_m$. Continuing in this fashion, we get $kh \in 
A_m$ for $k \in {\Bbb N}$. In particular, we reach the contradiction $0+e{\Bbb Z} = eh \in A_m$. 

This contradiction shows that, in fact, we have $H_m \cap B = \{0+e{\Bbb Z}\}$, so that $H_m+B 
\supseteq H_m \cup B$ and thus 
$$|H_m+B| \ge |H_m \cup B| = |H_m| + |B| - |H_m \cap B| = |H_m| + |B| -1.$$
Now Theorem 5 implies:
$$\eqalign{|A_m| &= |A_{m-1}+B| \ge |A_{m-1}+H_m| + |B+H_m| - |H_m| \cr
&\ge |A_{m-1}| + |B| -1 \ge (m-1)n+n = mn.\cr}$$
Thus $(\ast)$ holds. We assume that Theorem 1 is false. Then $m(q,e) > m$ where $m := \lceil {e  \over n} \rceil$.
This means that $0+e{\Bbb Z} \notin A_m$; in particular, we have $|A_m| \le |G|-1 = e-1$. On the other 
hand, $(\ast)$ implies that $|A_m| \ge mn \ge {e \over n}n = e$, and we have reached a contradiction.

%We note that, in the proof of Theorem, the precise nature of $G$ and $A_1$ is never used. Thus one may
%formulate more general versions of Theorem 1. 